\documentclass[reqno,a4paper,twoside]{article}
\usepackage{hyperref}
\usepackage{amsmath,amssymb,amsfonts,amsthm}

\newtheorem{The}{Theorem}[section]
\newtheorem{Pro}[The]{Proposition}

\newtheorem{Cor}[The]{Corollary}
\theoremstyle{definition}
\newtheorem{Def}[The]{Definition}
\newtheorem{Rem}[The]{Remark}
\newtheorem{Examp}[The]{Example}

\begin{document}

\title{Second order structures for sprays and connections on Fr\'{e}chet manifolds}

\author{M. Aghasi$^1$, A.R. Bahari$^1$, C.T.J. Dodson$^{2*}$,\\ 
G.N. Galanis$^3$ and A. Suri$^1$ \\
{\small $^1$Department of Mathematics,
Isfahan University of Technology,}\\ {\small Isfahan, Iran}\\
{\small $^2$School of Mathematics, University of Manchester,}\\ {\small Manchester M13 9PL, UK }\\
{\small $^3$Section of Mathematics, Naval Academy of Greece, 
Xatzikyriakion},\\ {\small Piraeus 185 39, Greece}}
\date{\small\it 29 October 2008}
\maketitle

\begin{abstract}Ambrose, Palais and Singer~\cite{Ambrose}
introduced the concept of second order structures on finite
dimensional manifolds. Kumar and Viswanath~\cite{Kumar}
extended these results to the category of Banach manifolds. In
the present paper all of these results are generalized to a large
class of Fr\'{e}chet manifolds. It is proved that the existence of
Christoffel and Hessian structures, connections, sprays and
dissections are equivalent on those Fr\'{e}chet manifolds which
can be considered as projective limits of Banach manifolds. These
concepts provide also an alternative way for the study of ordinary
differential equations on non-Banach infinite dimensional
manifolds. Concrete
examples of the  structures are provided using direct and flat connections.



$^*$Email: ctdodson@manchester.ac.uk
\end{abstract}

%%%%%%%%%%%%%%%%%%%%%%%%%%%%%%%%%%%%%%%%%%%%%%%%%%%%%%%%%%%%%%%    INtroduction  %%%%%%%%%%%%%%%%%%%%%%%%%%%%%%%%%%%%%%%%%%%%%%%%%%%%%%%%%%%%%%%%%%%

\section{introduction}The study of infinite dimensional manifolds
has received  much interest due to its interaction with
bundle structures, fibrations and foliations, jet fields, connections, sprays, Lagrangians
and Finsler structures (\cite{2de},\cite{EarleEells},\cite{Anto1}, \cite{Anto2},
\cite{Del}, \cite{Gal P} and \cite{Sau}). In particular,
non-Banach locally convex modelled manifolds have been studied from
different points of view (see for example \cite{conj},
\cite{Ab-Man}, \cite{Dod-Gal}, \cite{DGV05}, \cite{GalU} and
\cite{Neeb}). Fr\'{e}chet spaces of sections arise
naturally as configurations of a physical field and the
moduli space of inequivalent
configurations of a physical field is the quotient of the
infinite-dimensional configuration space $\mathcal{X}$ by the appropriate
symmetry gauge group. Typically, $\mathcal{X}$ is modelled on a Fr\'{e}chet
space of smooth sections of a vector bundle over a closed manifold.
For example, see Omori~\cite{Om,omori}.

The second order structures introduced by Ambrose et
al.~\cite{Ambrose} for finite dimensional manifolds were extended by Kumar and
Viswanath~\cite{Kumar} for Banach modelled manifolds.
They proved that Hessian structures, sprays,
dissections and (linear) connections are in a one-to-one
correspondence. However, there these concepts have to be supported by a
 Christoffel bundle and vector fields.
In this paper, following the lines of \cite{Kumar}, we first
construct the concepts of Christoffel bundle and fields for a
class of projective limit Fr\'{e}chet  manifolds. Then, we identify it with the other
structures, i.e. connections, Hessian structures and sprays.

One of the main problems in the study of non-Banach modelled
manifolds $M$ is the pathological structure of the general linear
group $GL(\mathbb{F})$ of a non-Banach space $\mathbb{F}.$ $GL(\mathbb{F})$
serves as the structure group of the tangent bundle $TM$, similar
to finite dimensional and Banach cases, but it is not even a
reasonable topological group structure
within the Fr\'{e}chet framework (see \cite{Gal 4}, \cite{Gal P}).\\
Moreover, for a Fr\'{e}chet space $\mathbb{F}$, $L(\mathbb{F}),$
the space of linear maps on $\mathbb{F},$ is not in general a Fr\'{e}chet
space. The same problem holds for the space of bilinear
maps $L^2(\mathbb{F},\mathbb{F})=\{B; B : \mathbb{F}
\times\mathbb{F} \longrightarrow \mathbb{F} , B $ is linear$\}$.

If one follows the classical procedure to define the notion of
Christoffel bundle or Hessian structures, then
$L^2(\mathbb{F},\mathbb{F})$ will appear as the corresponding
fibre type. As stated in Section~\ref{Chr}, these problems are
overcome by replacing $L^2(\mathbb{F},\mathbb{F})$ with an
appropriate Fr\'{e}chet space.
Another serious drawback in the study of Fr\'{e}chet manifolds and
bundles is the fact that there is no general solvability theory
for differential equations (\cite{Neeb}). This problem also can be
overcome if we restrict ourselves to the category of those
Fr\'{e}chet manifolds which can be considered as projective limits
of Banach corresponding factors.
 To eliminate these difficulties, we endow
$TM$ with a generalized vector bundle structure. (Note that
Galanis in \cite{Gal 4} proved a similar result but with a
different definition for tangent bundle). In the sequel we construct
the Christoffel bundles, connections, Hessian structures, sprays
and dissections. It is shown in this way that all the results
stated in \cite{Ambrose} and \cite{Kumar} hold in the
category of projective limit manifolds.

Our approach here gives the opportunity to study the problems
related to ordinary differential equations that arise via geometric
objects on manifolds. For example, geodesics with respect to
connections and sprays, and parallel transport are discussed.
Finally, the associated structures for flat and direct connections
are introduced.

%%%%%%%%%%%%%%%%%%%%%%%%%%%%%%%%%%%%%%%%%%%%%%%%%%%%%%%%%%%%%%%     Tangent bundle  %%%%%%%%%%%%%%%%%%%%%%%%%%%%%%%%%%%%%%%%%%%%%%%%%%%%%%%%%%%%%%%%%%%
%
%%%%%%%%%%%%%%%%%%%%%%%%%%%%%%%%%%%%%%%%%%%%%%%%%%%%%%%%%%%%%%%     Tangent bundle  %%%%%%%%%%%%%%%%%%%%%%%%%%%%%%%%%%%%%%%%%%%%%%%%%%%%%%%%%%%%%%%%%%%
%
%%%%%%%%%%%%%%%%%%%%%%%%%%%%%%%%%%%%%%%%%%%%%%%%%%%%%%%%%%%%%%%     Tangent bundle  %%%%%%%%%%%%%%%%%%%%%%%%%%%%%%%%%%%%%%%%%%%%%%%%%%%%%%%%%%%%%%%%%%%
\section{Christoffel bundle}\label{Chr}

Most of  our calculus is based on \cite{Abraham} and \cite{lang}.
Let $\mathbb{E}$ be a real Banach space, $M$ a Hausdorff
paracompact smooth manifold and  $m$ a point of $M$. The tangent
bundle of $M$ is defined as follows: $TM=\bigcup_{m\in M}T_mM$,
where $T_mM$ is considered as the set of equivalence classes of
all triples $(U,\varphi,e)$, where $(U,\varphi)$ is a chart of $M$
around $m$ and $e$ is an element of the model space $\mathbb{E}$
in which $\varphi U$ lies. $TM$ is a vector bundle on $M$ with
structure group $GL{({\mathbb{E}})}$ (\cite{lang}).

We summarise our basic notations about a certain rather wide class of
Fr\'{e}chet manifolds, namely those  which  can be considered as
projective limits of Banach manifolds. Let
$\{(M_i,\varphi_{ji})\}_{i,j\in\mathbb{N}}$ be a projective system
of Banach manifolds with $M=\varprojlim M_i$ such that for every $i
\in \mathbb{N}$, $M_i$ is modelled on the Banach space
$\mathbb{E}_i$ and $\{\mathbb{E}_i,\rho_{ji}\}_{i\in\mathbb{N}}$
 forms a projective system of Banach spaces. Furthermore
suppose that for each $m={(m)}_{i\in\mathbb{N}}\in M$ there exists
a projective system of local charts
$\{(U_i,\varphi_i)\}_{i\in\mathbb{N}}$ such that $m_i\in U_i$ and
$U=\varprojlim U_i$ is open in $M$ (see \cite{Ab-Man}).

%%%%%%%%%%%%%%%%%%%%%%%%%%%%%%%%%%%%%%%%%%%%%%%%%%%%%%%%%%%%%%%%%%%%%
It is known that for a Fr\'{e}chet space $\mathbb{F}$, the general
linear group $GL(\mathbb{F})$  cannot be endowed with a smooth
Lie group structure. It does not even admit a reasonable
topological group structure. The problems concerning the structure
group of $TM$ can be overcome by the replacement of
$GL(\mathbb{F})$ with the following topological group (and in a
generalized sense it is also a smooth Lie group):
\begin{eqnarray*}
\mathcal{H}_0{\mathbb{(F)}}=\{{(f_i)}_{i \in \mathbb{N}}\in
\prod_{i \in \mathbb{N}}GL(\mathbb{E}_i) :{ \varprojlim {f_i}} ~
\rm{exists} \}.
\end{eqnarray*}
More precisely $\mathcal{H}_0{\mathbb{(F)}}$ is isomorphic to the
projective limit of the Banach Lie groups
\begin{eqnarray*}
{\mathcal{H}_0}^i{\mathbb{(F)}}=\{(f_1,f_2,...,f_i) \in
\prod_{k=1}^{i }GL(\mathbb{E}_k): \rho_{jk}\circ
f_j=f_k\circ\rho_{jk}, {( k\leq j\leq i)}\}.
\end{eqnarray*}

%%%%%%%%%%%%%%%%%%%%%%%%%%%%%%%%%%%%%%%%%%%%%%%%%%%%%%%%%%%%
Under these notations the following basic theorems hold (compare
with \cite{Gal 4}).
\begin{The}
If $\{M_{i}\}_{i \in \mathbb{N}}$ is a projective system of
manifolds then $\{TM_i\}_{i \in \mathbb{N}}$ is also a projective
system with limit (set-theoretically) isomorphic to
$TM=\varprojlim TM_i$.
\end{The}

%%%%%%%%%%%%%%%%%%%%%%%%%%%%%%%%%%%%%%%%%%%%%%%%%%%%%%%%%%%%%%%%%%%%%

\begin{The}
$TM=\varprojlim TM_i$ has a Fr\'{e}chet vector bundle structure on
$M=\varprojlim M_i$ with structure group
$\mathcal{H}_0{\mathbb{(F)}}.$
\end{The}
%%%%%%%%%%%%%%%%%%%%%%%%%%%%%%%%%%%%%%%%%%%%%%%%%%%%%%%%%%%%%%%%
%%%%%%%%%%%%%%%%%%%%%%%%%%%%%%%%%%%%%%%%%%%%%%%%%%%%%%%%      The christoffel bundle          %%%%%%%%%%%%%%%%%%%%%%%%%%%%%%%%%%%%%%%%%%%%%%%%
%
%%%%%%%%%%%%%%%%%%%%%%%%%%%%%%%%%%%%%%%%%%%%%%%%%%%%%%%%       The christoffel bundle           %%%%%%%%%%%%%%%%%%%%%%%%%%%%%%%%%%%%%%%%%%%%%%%%
%
%%%%%%%%%%%%%%%%%%%%%%%%%%%%%%%%%%%%%%%%%%%%%%%%%%%%%%%%%      The christoffel bundle            %%%%%%%%%%%%%%%%%%%%%%%%%%%%%%%%%%%%%%%%%%%%%%%%
Let $L(\mathbb{E},\mathbb{E})$ be the space of continuous linear
maps from a Banach space $\mathbb{E}$ to  $\mathbb{E}$ and let
$L^2(\mathbb{E},\mathbb{E})$ be the space of all continuous
bilinear maps from $\mathbb{E}\times \mathbb{E}$ to $\mathbb{E}$.
For $m \in M$ and every chart $(U,\varphi)$ at $m$, consider the
triples of the form $(U,\varphi,B)$ where $B \in
L^2(\mathbb{E},\mathbb{E})$.

%%%%%%%%%%%%%%%%%%%%%%%%%%%%%%%%%%%%%%%%%%%%%%%%%%%%%%%%%%%%%%%%

\begin{Def}
Two triples $(U,\varphi,B_1)$ and $(V,\psi,B_2)$ are called
equivalent at $m$ if
\begin{eqnarray}
B_2{(DF(u).e_1, DF(u).e_2)}=DF(u)  . B_1(e_1,e_2) +
D^{2}F(u)(e_1,e_2),
\end{eqnarray}
where $u=\varphi m$, $F=\psi\circ {\varphi}^{-1}$ and $e_1,e_2 \in
\mathbb{E}$.
\end{Def}
%%%%%%%%%%%%%%%%%%%%%%%%%%%%%
It can  be checked that this is an equivalence relation. Each
equivalence class is called  a Christoffel element at $m$ and a
typical element is denoted by $\gamma$.  Let $(U,\varphi)$ be a
fixed chart at $m$. Define the mapping
\begin{eqnarray*}
C_\varphi :C_m \longmapsto L^2(\mathbb{E},\mathbb{E})\\
 \gamma\longmapsto {(\varphi{m},B)}
\end{eqnarray*}
where $C_m$ is the set of all Christoffel elements at $m$ and
$(U,\varphi,B)\in \gamma$. Then $C_\varphi$ is a bijection, which
endows  $CM=\bigsqcup_{m \in M} C_m$ with a $C^{\infty}$-atlas.
(For more details see \cite{Kumar}).

From \cite{Kumar} we have the result:
\begin{The}
The family $\{(CU,C\varphi)$: $(U,\varphi)$ \rm{is a chart on}
$M$\} is a $C^{\infty}$-atlas for $CM$.
\end{The}

%%%%%%%%%%%%%%%%%%%%%%%%%%%%%%%%%%%%%%%%%%%%%%%%%%%%%%%%%%%%%%%%%%%%%
We emphasise again at this point that for a
Fr\'{e}chet space $\mathbb{F}$, $L^2(\mathbb{F},\mathbb{F})$ does
not need to be a Fr\'{e}chet space in general. Hence,
the classical procedure for $CM$ for a
non-Banach Fr\'{e}chet manifold  $M$, does not yield a Fr\'{e}chet
manifold (nor bundle) structure. To overcome this obstacle we
use the Fr\'{e}chet space:
\begin{eqnarray*}
\mathcal{H}^2(\mathbb{F},\mathbb{F}):=\{(B_i)_{i \in \mathbb{N}}
\in \prod_{i \in \mathbb{N}} L^2(\mathbb{E}_i,\mathbb{E}_i):
~\varprojlim B_i~ \rm{exists} \}.
\end{eqnarray*}
$\mathcal{H}^2(\mathbb{F},\mathbb{F})$ is isomorphic to the
projective limit of Banach spaces
\begin{eqnarray*}
\mathcal{H}^2_i(\mathbb{F},\mathbb{F}):=\{(B_1,...,B_i) \in
\prod^i_{k=1}
L^2(\mathbb{E}_k,\mathbb{E}_k): ~B_k\circ(\rho_{jk}\times
\rho_{jk})=\rho_{jk}\circ B_j , {( k\leq j\leq i)}\}.
\end{eqnarray*}

Let $\{M_i\}_{i \in \mathbb{N}}$ be a projective system of Banach
manifolds as introduced earlier,  $B$, $\bar{B} \in
\mathcal{H}^2(\mathbb{F},\mathbb{F})$ and $(U=\varprojlim
U_i,\varphi=\varprojlim \varphi_i)$, $(V=\varprojlim V_i,
\psi=\varprojlim \psi_i)$ two corresponding charts.
\begin{Def}
Two triples $[U,\varphi,B]$ and $[V,\psi,\bar{B}]$ are equivalent
if, for every $i \in \mathbb{N}$, $[U_i,\varphi_i,B_i]$ and
$[V_i,\psi_i,{\bar{B}}_i]$ are equivalent.
\end{Def}
By these means one can show that $CM$ is endowed with a
Fr\'{e}chet manifold structure modelled on $\mathbb{F} \times
\mathcal{H}^2(\mathbb{F},\mathbb{F})$.
%%%%%%%%%%%%%%%%%%%%%%%%%%%%%%%%%%%%%%%%%%%%%%%%%%%%%%%%%%%%%%%%%%%%%%

\begin{Pro}
If $\{M_{i}\}_{i \in \mathbb{N}}$ is a projective system of
manifolds and  $\varprojlim CM_i$ exists then $\varprojlim
CM_i=C(\varprojlim M_i)$ (set-theoretically).
\begin{proof}
If we consider
\begin{eqnarray*}
Q:C(\varprojlim M_i)&\longrightarrow& \varprojlim (CM_i)\\
{[U,\varphi,B]} &\longmapsto &{ { ({[U_i,\varphi_i,B_i]}_{i})}_{i
\in \mathbb{N}}}
\end{eqnarray*}
then $Q$ is well defined. $Q$ is one to one since
$Q([U,\varphi,B])= Q([\bar{U},\bar{\varphi},\bar{B}])$ yields;
\begin{eqnarray*}
[U_i,\varphi_i,B_i]_i=[\bar{U_i},\bar{\varphi_i},\bar{B_i}]_i ~
,~i \in \mathbb{N}.
\end{eqnarray*}
Consequently $[U,\varphi,B]=[\varprojlim U_i,\varprojlim
\varphi_i,\varprojlim B_i]=\varprojlim[U_i,\varphi_i,B_i]_i=
\varprojlim[\bar{U_i},\bar{\varphi_i},\bar{B_i}]_i$\\
$=[\varprojlim \bar{U_i},\varprojlim \bar{\varphi_i},\varprojlim
\bar{B_i}]= [\bar{U},\bar{\varphi},\bar{B}]$.
Then $Q$ is also surjective since for every  $([U_i,\varphi_i,B_i]_i)_{i
\in \mathbb{ N}}$ in $\varprojlim (CM_i)$,
$Q(a)=([U_i,\varphi_i,B_i]_i)_{i \in \mathbb{ N}}$\\
 where
$a=[\varprojlim U_i,\varprojlim \varphi_i,\varprojlim B_i]$.

Therefore,  $Q$ is a bijection between $CM$ and $\varprojlim
(CM_i)$.
\end{proof}
\end{Pro}

%%%%%%%%%%%%%%%%%%%%%%%%%%%%%%%%%%%%%%%%%%%%%%%%%%%%%%%%%%%%%%%%%

The functions
\begin{eqnarray*}
\xi_\alpha:{\pi}^{-1}{(U_\alpha)}&\longrightarrow &U_\alpha\times
L^2(\mathbb{E},\mathbb{E})\\
\gamma &\longmapsto &{(m,B)};~ \alpha \in I
\end{eqnarray*}
with $\gamma \in C_m$, $(U_\alpha,\varphi_\alpha,B) \in \gamma$,
define a family of trivializations under which $(CM,M,\pi)$
becomes a fibre bundle ($\pi$ is the natural projection).

%%%%%%%%%%%%%%%%%%%%%%%%%%%%%%%%%%%%%%%%%%%%%%%%%%%%%%%%%%%%%%%%%%%%

In the next theorem the concept of $(CM,M,\pi)$ is generalized to a
Fr\'{e}chet manifold $M=\varprojlim M_i$.

\begin{The}
If $CM=\varprojlim CM_i$ exists, then it admits a Fr\'{e}chet
fibre bundle structure on $M=\varprojlim M_i$ with fibre type
$\mathcal{H}^2(\mathbb{F},\mathbb{F})$.
\begin{proof}
Let $\mathcal{A}=\{{(U_\alpha=\varprojlim
{U_\alpha}^{i},\varphi_\alpha=\varprojlim
{\varphi_\alpha}^{i})}\}$ be an atlas for $M=\varprojlim M_i$.
Then, for every $i \in \mathbb{N}$, $(CM_i,M_i,\pi_i)$ is a fibre
bundle with fibres of type $L^2(\mathbb{E}_i, \mathbb{E}_i)$ and
trivializations the mappings:
\begin{eqnarray*}
{\xi_\alpha}^{i}:{\pi}^{-1}_i{({U_\alpha}^{i})}&\longrightarrow
&{U_\alpha}^i\times
L^2(\mathbb{E}_i,\mathbb{E}_i)\\
\gamma_i &\longmapsto &{(m_i,B_i)}
\end{eqnarray*}
Suppose that $\{c_{ji}\}_{i,j \in \mathbb{N}}$,
$\{\varphi_{ji}\}_{i,j \in \mathbb{N}}$ and $\{\rho_{ji} \}_{i,j
\in \mathbb{N}}$  are the connecting morphisms of the projective
systems $CM=\varprojlim CM_i$, $M=\varprojlim M_i$ and
$\mathbb{F}=\varprojlim {\mathbb{E}_i}$  respectively. Since
$\varphi_{ji}\pi_j=\pi_i c_{ji}$, $\{\pi_i\}_{i \in \mathbb{N}}$
is a projective system of maps.  For every $\alpha \in {I}$,
$\{{\xi_\alpha}^i\}_{i \in \mathbb{N}}$ is a projective system and
$\pi=\varprojlim \pi_i:CM\longrightarrow M$ serves as the
projection map. On the other hand, $\xi_\alpha:=\varprojlim
{\xi_\alpha}^i:{\pi}^{-1}{({U_\alpha})}\longrightarrow {U_\alpha}
\times{\mathcal{H}^2(\mathbb{F},\mathbb{F})}$ is a diffeomorphism
since it is a projective limit of diffeomorphisms.
\end{proof}
\end{The}

%%%%%%%%%%%%%%%%%%%%%%%%%%%%%%%%%%%%%%%%%%%%%%%%%%%%%%%%%%%%%%%

For an  open subset $U$ in $\mathbb E$, define a Christoffel map
$\Gamma$ on $U$ to be a smooth mapping $\Gamma:U\longrightarrow
L^2(\mathbb{E}, \mathbb{E})$ and for every chart $(U,\varphi)$ of
$M$ a Christoffel map is locally a smooth mapping
 $\Gamma_\varphi:\varphi U\longrightarrow
L^2(\mathbb{E}, \mathbb{E}).$

%%%%%%%%%%%%%%%%%%%%%%%%%%%%%%%%%%%%%%%%%%%%%%%%%%%%%%%%%%

\begin{Def}
$M$ is endowed with a Christoffel structure $\{\Gamma_\varphi\}$
if for every choice of charts $(U,\varphi)$ and $(V,\psi)$ with
$U\cap V\neq\O$, the following relation holds true:
\begin{eqnarray*}
\Gamma_\psi{(v)}{(DF{(u)}.e_1,DF{(u)}.e_2)}=DF{(u)}.\Gamma_\varphi{(u)}
{(e_1,e_2)}+D^2F{(u)}.{(e_1,e_2)}
\end{eqnarray*}
where $e_1,e_2\in \mathbb E$, $\varphi m=u$, $\psi m=v$ and
$F=\psi\circ {\varphi}^{-1}$.

For $v,w\in T_mM$ we can express this condition as follows:
\begin{eqnarray*}
\Gamma_\psi{(v)}{(v_\psi,w_\psi)}=DF{(u)}.\Gamma_\varphi{(u)}
{(v_\varphi,w_\varphi)}+D^2F{(u)}.{(v_\varphi,w_\varphi)}
\end{eqnarray*}
where $\varphi m=u$, $v_\psi=DF(u) . v_\varphi$, $w_\psi=DF(u) .
w_\varphi$, $v=[U,\varphi,v_\varphi]$ and
$w=[U,\varphi,w_\varphi]$ (see also \cite{Kumar}).
\end{Def}

%%%%%%%%%%%%%%%%%%%%%%%%%%%%%%%%%%%%

In a similar manner  one can define the Christoffel map for the
non-Banach case as follows: Let $U=\varprojlim U_i$ be an open
subset of $\mathbb {F}=\varprojlim
 \mathbb{E}_i$. A Christoffel map  on $U=\varprojlim U_i$, is a projective limit smooth mapping $\Gamma=\varprojlim \Gamma_i:U\longrightarrow
{\mathcal{H}^2(\mathbb{F},\mathbb{F})} $. Note that for each chart
$(U=\varprojlim U_i,\varphi=\varprojlim \varphi_i)$ of $M$,
$\varprojlim \Gamma_{\varphi_i}:=\Gamma_\varphi:\varphi
U\longrightarrow {\mathcal{H}^2(\mathbb{F},\mathbb{F})}$ defines a
Christoffel map on $U$. Now we can state the following definition
for Fr\'{e}chet manifolds.

%%%%%%%%%%%%%%%%%%%%%%%%%%%%%%%%%%

\begin{Def}
$M=\varprojlim M_i$ is endowed with a Christoffel structure
$\{\Gamma_\varphi=\varprojlim \Gamma_{\varphi_i}\}$, if for every pair of  charts $(U=\varprojlim U_i,\varphi=\varprojlim
\varphi_i)$ and $(V=\varprojlim V_i,\psi=\varprojlim {\psi_i})$
around $m=(m_i)_{i \in \mathbb{N}}$ the following relation is
satisfied:
\begin{eqnarray*}
\Gamma_\psi{(v)}{(DF{(u)}.e_1,DF{(u)}.e_2)}=DF{(u)}.\Gamma_\varphi{(u)}
{(e_1,e_2)}+D^2F{(u)}.{(e_1,e_2)},
\end{eqnarray*}
where $e_1=(e^i_1)_{i \in \mathbb{N}}, e_2=(e^i_2)_{i \in
\mathbb{N}} \in \mathbb{F}$, $\varprojlim \varphi_i
m_i=\varprojlim u_i=u$, $\varprojlim \psi_i m_i=\varprojlim v_i=v$
and $F=\varprojlim F_i=\varprojlim \psi_i\circ {\varphi_i}^{-1}$.
For $v,w\in T_mM=\varprojlim T_{m_i}M_i$ this condition takes the
form
\begin{eqnarray*}
\Gamma_\psi{(v)}{(v_\psi,w_\psi)}=DF{(u)}.\Gamma_\varphi{(u)}
{(v_\varphi,w_\varphi)}+D^2F{(u)}.{(v_\varphi,w_\varphi)}
\end{eqnarray*}
where $\varprojlim \varphi_i m_i=\varprojlim u_i=u$,
$v_\psi=\varprojlim DF_i(u_i) . v_{\varphi_i}$,
$w_\psi=\varprojlim DF_i(u_i) . w_{\varphi_i}$,
$v=([U_i,\varphi_i,v_{\varphi_i}]_i)_{i \in \mathbb{N}}$ and
$w=([U_i,\varphi_i,w_{\varphi_i}]_i)_{i \in\mathbb{N}}$.
\end{Def}

%%%%%%%%%%%%%%%%%%%%%%%%%%%%%%%%%%%%%%%%%%%%%%%%%%%%%%%%       Connections             %%%%%%%%%%%%%%%%%%%%%%%%%%%%%%%%%%%%%%%%%%%%%%%%
%
%%%%%%%%%%%%%%%%%%%%%%%%%%%%%%%%%%%%%%%%%%%%%%%%%%%%%%%%       Connections             %%%%%%%%%%%%%%%%%%%%%%%%%%%%%%%%%%%%%%%%%%%%%%%%
%
%%%%%%%%%%%%%%%%%%%%%%%%%%%%%%%%%%%%%%%%%%%%%%%%%%%%%%%%       Connections             %%%%%%%%%%%%%%%%%%%%%%%%%%%%%%%%%%%%%%%%%%%%%%%%
\section{Connections and Hessian structures}
%%%%%%%%%%%%%%%%%%%%%%%%%%%%
A connection on $M$ by Koszul's definition (see \cite{Flaschel})
is a smooth mapping
\begin{eqnarray*}
\nabla :\chi(M)\times\chi(M)&\longrightarrow & \chi(M)\\
(X,Y) &\longmapsto& \nabla_XY
\end{eqnarray*}
such that on every local chart $(U,\varphi)$ on $M$, there exists a
smooth map $\Gamma_\varphi:\varphi U\longrightarrow
L^2(\mathbb{E}, \mathbb{E})$ with
\begin{eqnarray*}
(\nabla_XY)(\varphi m)=DY_\varphi(\varphi m).X_\varphi(\varphi m)
-\Gamma_\varphi(\varphi m)(X_\varphi(\varphi m),Y_\varphi(\varphi
m));~ \forall m \in U.
\end{eqnarray*}

%%%%%%%%%%%%%%%%%%%%%%%%%%%%%%%%

We prove in the sequel that if $\nabla$ is a connection on
$M$, then $\{\Gamma_\varphi\}$ forms a Christoffel structure on
$M$. Conversely if $\{\Gamma_\varphi\}$ is a Christoffel structure
on $M$ and $X,Y \in \chi(U)$, then a connection $\nabla$ can be
defined by
\begin{eqnarray*}
(\nabla_XY)(m)=T{\varphi}^{-1}[DY_\varphi(\varphi
m).X_\varphi(\varphi m) -\Gamma_\varphi(\varphi
m)(X_\varphi(\varphi m),Y_\varphi(\varphi m))]
\end{eqnarray*}
(see \cite{Kumar}).

%%%%%%%%%%%%%%%%%%%%%%%%%%%%%%%

Before proceeding to results, it is necessary to prove the
following.
\begin{The}\label{projcon}
The limit $\nabla=\varprojlim \nabla_i$ of a projective system of
connections $\{\nabla_i\}_{i \in \mathbb{N}}$  is a connection on
$M=\varprojlim M_i$.
\begin{proof}
For $i\leq j$, let $(U_j,\varphi_j)$ be a chart of $M_j$ around
$m_j$ and $(U_i,\varphi_i)$ be a chart of $M_i$ at $\varphi_{ji}
{m_j}=m_i$. Moreover for every $i \in \mathbb{N}$, let
$X_{\varphi_i}:\varphi_i{U_i}\longrightarrow \mathbb{E}_i$ be the
local principal part of  $X_i \in \chi(M_i)$. Since $\nabla$ is a
smooth mapping as a projective limit of smooth factors, to prove
the theorem it suffices to check that
%%%%%%%%%%%%%%%%%%
${\rho_{ji}}\circ \nabla_{X_{\varphi_j}}Y_{\varphi_j}
=\nabla_{X_{\varphi_i}}Y_{\varphi_i}\circ{\rho_{ji}}$.\\
%%%%%%%%%%%%%%%%%%
The last equality holds since for $m_j\in U_j$;
\begin{eqnarray*}
&&{\rho_{ji}}\circ \nabla_{X_{\varphi_j}}Y_{\varphi_j}(\varphi_j{m_j})\\
&=&\underbrace{{\rho_{ji}}DY_{\varphi_j}(\varphi_j{m_j}).X_{\varphi_j}(\varphi_j{m_j})}_{*}
-\underbrace{\rho_{ji}\Gamma_{\varphi_j}(\varphi_j{m_j})
(X_{\varphi_j}(\varphi_j{m_j}),Y_{\varphi_j}(\varphi_j{m_j}))}_{**}\\
&=&\nabla_{X_{\varphi_i}}Y_{\varphi_i}(\varphi_i{m_i})=
\nabla_{X_{\varphi_i}}Y_{\varphi_i}\rho_{ji}(\varphi_j{m_j})
\end{eqnarray*}

Note that
\begin{eqnarray*}
*&=&\frac d{dt}\rho_{ji}Y_{\varphi_j}(\varphi_j{m_j} +
tX_{\varphi_j}(\varphi_j{m_j}))_{|t=0}\\
%%%%%%%%%%%%%%%%%%%%%%%%%%%%%
&=&\frac d{dt}Y_{\varphi_i}\rho_{ji}(\varphi_j{m_j} +
tX_{\varphi_j}(\varphi_j{m_j})))_{|t=0}\\
%%%%%%%%%%%%%%%%%%%%%%%%%%%%%
&=&\frac d{dt}Y_{\varphi_i}(\varphi_i{m_i} +
tX_{\varphi_i}(\varphi_i{m_i}))_{|t=0}\\
%%%%%%%%%%%%%%%%%%%%%%%%%%%%%
&=&DY_{\varphi_i}(\varphi_i{m_i)}.X_{\varphi_i}(\varphi_i{m_i})
\end{eqnarray*}

and
\begin{eqnarray*}
**&=&\Gamma_{\varphi_i}(\varphi_i{m_i})(\rho_{ji}\times\rho_{ji})
(X_{\varphi_j}(\varphi_j{m_j}),Y_{\varphi_j}(\varphi_j{m_j}))\\
&=&\Gamma_{\varphi_i}(\varphi_i{m_i})(X_{\varphi_i}
(\varphi_i{m_i}),Y_{\varphi_i}(\varphi_i{m_i}))
\end{eqnarray*}
\end{proof}
\end{The}

%%%%%%%%%%%%%%%%%%%%%%%%%%%%%%%%%%%%%%%%%%%%%%%%%%%%%%%%%%%%%%%%%%%%%%

Based on Theorem~\ref{projcon} we may now establish several important
properties.
\begin{The}\label{ChrStr}
If $\nabla=\varprojlim \nabla_i$ is a connection on $M=\varprojlim
M_i$, then $\{\Gamma_{\varphi}=\varprojlim \Gamma_{\varphi_i}\}$
forms a Christoffel structure on $M$.
\begin{proof}
Let $(U$=$\varprojlim U_i$, $\varphi$=$\varprojlim \varphi_i)$,
$(V$=$\varprojlim v_i$, $\psi$=$\varprojlim \psi_i)$ be two charts
through  $m=(m_i)_{i \in \mathbb{N}} \in M$ and  $\varprojlim
\varphi_i{m_i}=\varprojlim u_i=u$, $F=\varprojlim F_i=\varprojlim
(\psi_i\circ {\varphi_i}^{-1})$. Furthermore, suppose that
$X_\varphi=\varprojlim X_{\varphi_i}$, then
\begin{eqnarray*}
&&[DY_\psi.X_\psi](F(u))
%%%%%%%%%%%%%%%%%%%
=[D\varprojlim Y_{\psi_i}.\varprojlim X_{\psi_i}](\varprojlim
F_i(\varprojlim u_i))
%%%%%%%%%%%%%%%%%%%%%
=\varprojlim[[DY_{\psi_i}.X_{\psi_i}](F_i(u_i))]\\
%%%%%%%%%%%%%%%%%%%%%%
&=&\varprojlim [DY_{\psi_i}(F_i(u_i)).X_{\psi_i}(F_i(u_i))]
%%%%%%%%%%%%%%%%%%%%%%
=\varprojlim [DY_{\psi_i}(F_i(u_i)).DF_i(u_i).X_{\varphi_i}(u_i)]\\
%%%%%%%%%%%%%%%%%%%%%
&=&\varprojlim [D(Y_{\psi_i}\circ F_i)(u_i).X_{\varphi_i}(u_i)]
%%%%%%%%%%%%%%
=\varprojlim [D(DF_i.Y_{\varphi _i}(u_i).X_{\varphi_i}(u_i)]\\
%%%%%%%%%%%%%%
&=&\varprojlim
[D^2F_i(u_i)(X_{\varphi_i},Y_{\varphi_i})+DF_i(u_i).DY_{\varphi_i}(u_i).X_{\varphi_i}(u_i)].
\end{eqnarray*}
But
\begin{eqnarray*}
(\nabla_XY)_\varphi\circ F=\varprojlim
[(\nabla_{X_i}{Y_i})_{\varphi_i}\circ F_i]
%%%%%%%%%%%%%%%%%%%%%%%%
=\varprojlim [(DY_{\psi_i}.X_{\psi_i})\circ
F_i-\Gamma_{\psi_i}(X_{\psi_i},Y_{\psi_i})]\\
%%%%%%%%%%%%%%%%%%%%%%%%
=\varprojlim
[D^2F_i(X_{\varphi_i},Y_{\varphi_i})+DF_i.(DY_{\varphi_i}.X_{\varphi_i})
-\Gamma_{\psi_i}(X_{\psi_i},Y_{\psi_i})],
%%%%%%%%%%%%%%%%%%%%%%%%%
\end{eqnarray*}
hence
\begin{eqnarray*}
\Gamma_{\psi}(X_{\psi},Y_{\psi})=\varprojlim
[\Gamma_{\psi_i}(X_{\psi_i},Y_{\psi_i})]=
%%%%%%%%%%%%%%%%%%%%%%%%%
\varprojlim
[D^2F_i(X_{\varphi_i},Y_{\varphi_i})+\\DF_i.(DY_{\varphi_i}.X_{\varphi_i})
%%%%%%%%%%%%%%%%%%%%%%%%%
-DF_i.(\nabla_{X_i}{Y_i})_{\varphi_i}]
%%%%%%%%%%%%%%%%
=D^2F(X_\varphi,Y_\varphi)+DF.\Gamma_{\varphi}(X_\varphi,Y_\varphi).
\end{eqnarray*}
i.e. $\{\Gamma_\varphi=\varprojlim \Gamma_{\varphi_i}\}$
 forms a Christoffel structure on $M=\varprojlim M_i$.
\end{proof}
\end{The}

%%%%%%%%%%%%%%%%%%%%%%%%%%%%%%%%%%%%%%%%%%%%%%%%%%%%%%%%%%%%%%

\begin{Rem}
The converse also of Theorem~\ref{ChrStr} can be obtained by setting
\begin{eqnarray*}
(\nabla _XY)(m)=\varprojlim [T{\varphi^{-1}
_i}[DY_{\varphi_i}(\varphi_i{m_i}).X_{\varphi_i}(\varphi_i{m_i})-\\\Gamma_{\varphi_i}
(\varphi_i{m_i})(X_{\varphi_i}(\varphi_i{m_i}),Y_{\varphi_i}(\varphi_i{m_i}))]],
\end{eqnarray*}
where $\{\Gamma_{\varphi}=\varprojlim \Gamma_{\varphi_{i}}\}$ is a
Christoffel structure on $M$. Moreover for $f=\varprojlim f_i \in
C^{\infty}(M)$ and $X,Y \in \chi(M), \nabla$ satisfies the
following conditions:
\\
$(i)\nabla$ is real linear in $X$ and $Y$,\\
$(ii)\nabla_{fX}Y=f\nabla_XY$,\\
$(iii)\nabla_X(fY)=f\nabla_XY+(Xf)Y$.
\end{Rem}
%%%%%%%%%%%%%%%%%%%%%%%%%%%%%%%%%%%%%%%%%%%%%%%%%%%%%%%%%%%%%%%%% Hessian structures %%%%%%%%%%%%%%%%%%%%%%%%%%%%%%%%%%%%%%%%%%%%%%%%%%%%%%%%%%%%%%%
%
%%%%%%%%%%%%%%%%%%%%%%%%%%%%%%%%%%%%%%%%%%%%%%%%%%%%%%%%%%%%%%%%% Hessian structures %%%%%%%%%%%%%%%%%%%%%%%%%%%%%%%%%%%%%%%%%%%%%%%%%%%%%%%%%%%%%%%
%
%%%%%%%%%%%%%%%%%%%%%%%%%%%%%%%%%%%%%%%%%%%%%%%%%%%%%%%%%%%%%%%%% Hessian structures %%%%%%%%%%%%%%%%%%%%%%%%%%%%%%%%%%%%%%%%%%%%%%%%%%%%%%%%%%%%%%%
%%%%%%%%%%%%%%%%%%%%%%%%%%%%%%%%%%%%%%%%%%%%%%%%%%%%%%%%%%%%%%%%%%%%%%%%
In anticipation of the sequel, a Hessian structure on $M$ is a
mapping $H : f \longmapsto Hf$, which associates to every $f \in
C^{\infty}(M)$ a covariant 2-tensor $Hf$ on $M$ such that on a
local chart $(U,\varphi)$ of $M$ and for every $X,Y \in \chi(M)$,
there exists a smooth map $\Gamma_\varphi:\varphi U\longrightarrow
L^2(\mathbb{E},\mathbb{E})$ with
\begin{eqnarray*}
[Hf(X,Y)]_\varphi(\varphi m)
%%%%%%%%%%%%%%%%%%%%
 =D^2f_\varphi(\varphi m)(X_\varphi
(\varphi m),Y_\varphi (\varphi m))+\\
%%%%%%%%%%%%%%%%%%%%
Df_\varphi(\varphi m).\Gamma_\varphi(\varphi m)(X_\varphi (\varphi
m),Y_\varphi (\varphi m)).
\end{eqnarray*}

%%%%%%%%%%%%%%%%%%%

It turns out that $Hf$ is a Hessian structure on $M$ if and only
if $M$ admits the Christoffel structure $\{\Gamma_\varphi\}$.
Moreover, there is a one-to-one correspondence between Hessian
structures and connections given by
$Hf(X,Y)=X(Y(f))-(\nabla_{X}Y)f$. (For more details see
\cite{Kumar}).

%%%%%%%%%%%%%%%%%%%%%%%%%%%

Here we study the above results for projective limit
manifolds. However, we should consider just the smooth functions and
smooth vector fields such that $\mathcal{F}(M)=\{(f_i)_{i \in
\mathbb{N}}: f_i:M_i \longrightarrow \mathbb{R}$ is continuous and
$\varprojlim f_i$ exists$\}$ and $\mathcal{G}(M)=\{(X_i)_{i
\in\mathbb{N}}: X_i$ is a vector field on $M_i$ and $\varprojlim
X_i$ exists$\}$ respectively.

\begin{Pro}
The limit of a projective system of Hessian structures on
$\{M_i\}_{i \in \mathbb{N}}$ is a Hessian structure on
$M=\varprojlim M_i$.
\begin{proof}
 For every $i \in \mathbb{N}$, let $f_i \in C^\infty(M_i)$ and $X_i,Y_i \in
\chi(M_i)$. Consider a chart $(U_i,\varphi_i)$ on $M_i$. Assume
that $\Gamma_{\varphi_i}:\varphi_i{U_i}\longrightarrow
L^2(\mathbb{E}_i,\mathbb{E}_i)$ is a smooth map such that
\begin{eqnarray*}
[H_if_i(X_i,Y_i)]_{\varphi_i}(\varphi_i{m_i})
%%%%%%%%%%%%%%%%
=D^2{f_i}_{\varphi_i}(\varphi_i{m_i})({X}_{\varphi_i}
(\varphi_i{m_i}),{Y}_{\varphi_i} (\varphi_i{m_i}))+\\
%%%%%%%%%%%%%%%%
D{f_i}_{\varphi_i}({\varphi_i}{m_i}).{\Gamma}_{\varphi_i}
({\varphi_i}{m_i})({X}_{\varphi_i}
(\varphi_i{m_i}),{Y}_{\varphi_i} (\varphi_i{m_i})).
\end{eqnarray*}

Hence we must check that for $j\geq i$,
$[H_jf_j(X_j,Y_j)]_{\varphi_j}=
[H_if_i(X_i,Y_i)]_{\varphi_i}\circ\rho_{ji}$. For $m_j \in U_j$;
\begin{eqnarray*}
[H_jf_j(X_j,Y_j)]_{\varphi_j}(\varphi_j{m_j})
%%%%%%%%%%%%%%%%%%
=\underbrace{D^2{f_j}_{\varphi_j}(\varphi_j{m_j})({X}_{\varphi_j}
(\varphi_j{m_j}),{Y}_{\varphi_j} (\varphi_j{m_j}))}_{*}+\\
%%%%%%%%%%%%%%%%%%
\underbrace{D{f_j}_{\varphi_j}({\varphi_j}{m_j})
.{\Gamma}_{\varphi_j}({\varphi_j}{m_j})({X}_{\varphi_j}
(\varphi_j{m_j}),{Y}_{\varphi_j} (\varphi_j{m_j}))}_{**}.\\
%%%%%%%%%%%%%%%%%%
=[H_if_i(X_i,Y_i)]_{\varphi_i}(\varphi_i{m_i})=
%%%%%%%%%%%%%%%%%%
[H_if_i(X_i,Y_i)]_{\varphi_i}{\rho_{ji}}(\varphi_j{m_j}).
\end{eqnarray*}
Note that;
\begin{eqnarray*}
&&D{f_j}_{\varphi_j}(\varphi_j{m_j})({X}_{\varphi_j}
(\varphi_j{m_j}))
%%%%%%%%%%%%%%%%
=D({f_i}_{\varphi_i}\circ\rho_{ji})(\varphi_j{m_j})({X}_{\varphi_j}
(\varphi_j{m_j}))\\
%%%%%%%%%%%%%%%%
&=&{\frac
d{dt}}({f_i}_{\varphi_i}\circ\rho_{ji})(\varphi_j{m_j}+t{X}_{\varphi_j}
(\varphi_j{m_j}))|_{t=0}\\
%%%%%%%%%%%%%%%%
&=&{\frac d{dt}}{f_i}_{\varphi_i}(\varphi_i{m_i}+
t{X}_{\varphi_i}\rho_{ji}(\varphi_j{m_j}))|_{t=0}\\
%%%%%%%%%%%%%%%
&=&{\frac d{dt}}{f_i}_{\varphi_i}(\varphi_i{m_i}+ t{X}_{\varphi_i}
(\varphi_i{m_i}))|_{t=0}
%%%%%%%%%%%%%%%
=D{f_i}_{\varphi_i}(\varphi_i{m_i})({X}_{\varphi_i}
(\varphi_i{m_i})),
%%%%%%%%%%%%
\end{eqnarray*}
and consequently
\begin{eqnarray*}
*&=&D^2(f_j\circ{\varphi_j}^{-1})(\varphi_j{m_j})({X}_{\varphi_j}
(\varphi_j{m_j}),{Y}_{\varphi_j} (\varphi_j{m_j}))\\
%%%%%%%%%%%%%%%%
&=&D\big{(}D({f_i}_{\varphi_i}\circ\rho_{ji})(\varphi_j{m_j})({X}_{\varphi_j}
(\varphi_j{m_j}))\big{)}[{Y}_{\varphi_j} (\varphi_j{m_j})]\\
%%%%%%%%%%%%%%%%
%%%%%%%%%%%%%%%
&=&D\big{(}D{f_i}_{\varphi_i}(\varphi_i{m_i})({X}_{\varphi_i}
(\varphi_i{m_i})\big{(}[{Y}_{\varphi_i} (\varphi_i{m_i})]\\
%%%%%%%%%%%%
&=&D^2{f_i}_{\varphi_i}(\varphi_i{m_i})({X}_{\varphi_i}
(\varphi_i{m_i}),{Y}_{\varphi_i} (\varphi_i{m_i})).
\end{eqnarray*}
Moreover
\begin{eqnarray*}
**&=&D({f_i}_{\varphi_i}\circ\rho_{ji})(\varphi_j{m_j})
.{\Gamma}_{\varphi_j}({\varphi_j}{m_j})({X}_{\varphi_j}
(\varphi_j{m_j}),{Y}_{\varphi_j} (\varphi_j{m_j}))\\
%%%%%%%%%%%%%%%%%%%%
&=&{\frac d{dt}}({f_i}_{\varphi_i}\circ\rho_{ji})(\varphi_j{m_j}
+t{\Gamma}_{\varphi_j}({\varphi_j}{m_j})({X}_{\varphi_j}
(\varphi_j{m_j}),{Y}_{\varphi_j} (\varphi_j{m_j}))|_{t=0}\\
%%%%%%%%%%%%%%%%%%%
&=&{\frac
d{dt}}{f_i}_{\varphi_i}(\varphi_i{m_i}+t{\Gamma}_{\varphi_i}({\varphi_i}{m_i})
(\rho_{ji}\times\rho_{ji})({X}_{\varphi_j}
(\varphi_j{m_j}),{Y}_{\varphi_j} (\varphi_j{m_j}))|_{t=0}\\
%%%%%%%%%%%%%%%%%%
&=&D{f_i}_{\varphi_i}({\varphi_i}{m_i})
.{\Gamma}_{\varphi_i}({\varphi_i}({m_i})({X}_{\varphi_i}
(\varphi_i({m_i}),{Y}_{\varphi_i} (\varphi_i{m_i})).
\end{eqnarray*}

Hence $\varprojlim
[H_if_i(X_i,Y_i)]_{\varphi_i}=[Hf(X,Y)]_{\varphi}$ where $f \in
\mathcal{F}(M)$, $\varphi=\varprojlim \varphi_i$ and $X,Y \in
\mathcal{G}(M)$.
\end{proof}
\end{Pro}

%%%%%%%%%%%%%%%%%%%%%%%%%%%%%%%%%%%%%%%%%%%%%%%%%%%%%%%%%%%%%%%%%%%

Next, Theorem~\ref{HeqCon} proves that there is a one-to-one
correspondence between Hessian structures and connections on
Fr\'{e}chet manifolds.

%%%%%%%%%%%%%%%%%%%%%%%%%%%%%%%%%%%%%%%%%%%%%%%%%%%%%%%%%%%%%%%%%%%%

\begin{The}\label{HeqCon}
Let $\nabla=\varprojlim \nabla_i$ be a connection on
$M=\varprojlim M_i$, and $Hf(X,Y):=X(Y(f))-(\nabla_XY)f$. Then $H$
is a Hessian structure on $M$. Conversely the connection which
obtained as projective limit of connections  arises from a Hessian
structure.
\begin{proof}
Let $v,w \in T_mM=\varprojlim T_{m_i}M_i$ and $(U$=$\varprojlim
U_{i}$,$\varphi$=$\varprojlim \varphi_{i}$ be a chart around
$m=(m_i)_{i \in \mathbb{N}}$. Consider vector fields $\varprojlim
X_i,\varprojlim Y_i \in \chi(\varprojlim U_{i})$ with $\varprojlim
X_i(m_i)=v$ and $\varprojlim Y_i(m_i)=w$. Suppose
$\nabla=\varprojlim \nabla_i$ be a connection on $M=\varprojlim
M_i$, then $\{\Gamma_\varphi=\varprojlim
\Gamma_{\varphi_i}\}$ is a  Christoffel structure on $M$. Hence\\
%\begin{eqnarray*}

 $X(Y(f))(m)-(\nabla_XY).
f(m)=X_\varphi(Y_\varphi(f_\varphi))(\varphi m)-(\nabla_{X}Y)_\varphi(f_\varphi)(\varphi m)\\
=\varprojlim
[X_{\varphi_i}(Y_{\varphi_i}(f_{\varphi_i}))({\varphi_i}{m_i})
-(\nabla_{X_i}{Y_i})_{\varphi_i}(f_{i_{\varphi_i}})(\varphi_i{m_i})]$\\
%%%%%%%%%%%%%%%%%%%%%%%
$=\varprojlim
[D(Df_{i_{\varphi_i}}.Y_{\varphi_i})(\varphi_i{m_i}).X_{\varphi_i}(\varphi_i{m_i})
-Df_{i_{\varphi_i}}(\varphi_i{m_i}).(\nabla_{X_i}{Y_i})_{\varphi_i}(\varphi_i{m_i})]\\
%%%%%%%%%%%%%%%%%%%%%%
=\varprojlim [D^2{f_i}_{\varphi_i}({\varphi_i}{m_i})
(X_{\varphi_i}{(\varphi_i}{m_i}),(Y_{\varphi_i}{(\varphi_i}{m_i})\\
%%%%%%%%%%%%%%%%%%%%%
+D{f_i}_{\varphi_i}({\varphi_i}{m_i}).
DY_{\varphi_i}{(\varphi_i}{m_i}).X_{\varphi_i}{(\varphi_i}{m_i})
%%%%%%%%%%%%%%%%%%%%%
-D{f_i}_{\varphi_i}({\varphi_i}{m_i})[DY_{\varphi_i}{(\varphi_i}{m_i}).X_{\varphi_i}{(\varphi_i}{m_i})\\
-\Gamma_{\varphi_i}({\varphi_i}{m_i})(X_{\varphi_i}{(\varphi_i}{m_i}),(Y_{\varphi_i}{(\varphi_i}{m_i})]]$\\
%\end{eqnarray*}
$=\varprojlim
[D^2{f_i}_{\varphi_i}(\varphi_i{m_i})({X}_{\varphi_i}
(\varphi_i{m_i}),{Y}_{\varphi_i} (\varphi_i{m_i}))+\\
%%%%%%%%%%%%%%%%
D{f_i}_{\varphi_i}({\varphi_i}{m_i}).{\Gamma}_{\varphi_i}
({\varphi_i}{m_i})({X}_{\varphi_i}
(\varphi_i{m_i}),{Y}_{\varphi_i} (\varphi_i{m_i}))]$\\
 $=\varprojlim
[[H_if_i(X_i,Y_i)]_{\varphi_i}(\varphi_i{m_i})]$
$=[Hf(X,Y)]_\varphi(\varphi{m})$.  \\
Conversely if $Hf=\varprojlim H_if_i$ is a Hessian structure on
$M=\varprojlim M_i$ then $\{\Gamma_{\varphi}=\varprojlim
\Gamma_{\varphi_i}\}$ forms a Christoffel structure on $M$. Now we have \\

%\begin{eqnarray*}
%%%%%%%%%%%%%%%%%%%%%%%%%%%%%%%%%
$X(Y(f))(m)-Hf(X,Y)(m)=X_\varphi(Y_\varphi(f_\varphi))(\varphi m)-
[Hf(X,Y)]_\varphi(\varphi{m})$\\
%%%%%%%%%%%%%%%%%%%%%%%%%%%%%%%%
$=\varprojlim
[X_{\varphi_i}(Y_{\varphi_i}(f_{\varphi_i}))({\varphi_i}{m_i})
-[H_if_i(X_i,Y_i)]_{\varphi_i}(\varphi_i{m_i})]$\\
%%%%%%%%%%%%%%%%%%%%%%%%%%%%%%%%
$=\varprojlim [D^2{f_i}_{\varphi_i}({\varphi_i}{m_i})
(X_{\varphi_i}{(\varphi_i}{m_i}),(Y_{\varphi_i}{(\varphi_i}{m_i})+
D{f_i}_{\varphi_i}({\varphi_i}{m_i}).
DY_{\varphi_i}{(\varphi_i}{m_i}).X_{\varphi_i}{(\varphi_i}{m_i})$\\
%%%%%%%%%%%%%%%%%%%%%%%%%%%%%%%%
$-D^2{f_i}_{\varphi_i}(\varphi_i{m_i})({X}_{\varphi_i}
(\varphi_i{m_i}),{Y}_{\varphi_i} (\varphi_i{m_i}))\\
%%%%%%%%%%%%%%%%
-D{f_i}_{\varphi_i}({\varphi_i}{m_i}).{\Gamma}_{\varphi_i}
({\varphi_i}{m_i})({X}_{\varphi_i}
(\varphi_i{m_i}),{Y}_{\varphi_i} (\varphi_i{m_i}))]$\\
$=\varprojlim
[Df_{i_{\varphi_i}}(\varphi_i{m_i}).(\nabla_{X_i}{Y_i})
_{\varphi_i}(\varphi_i{m_i})]=\varprojlim[(\nabla_{X_i}{Y_i})_{\varphi_i}
(f_{i_{\varphi_i}})(\varphi_i{m_i})]$\\
$=(\nabla_{X}Y)_\varphi(f_\varphi)(\varphi m)$
\end{proof}
\end{The}

%%%%%%%%%%%%%%%%%%%%%%%%%%%%%%%%%%%%%%%%%%%%%%%%%%%%%%%%%%%%%%%%% Spray %%%%%%%%%%%%%%%%%%%%%%%%%%%%%%%%%%%%%%%%%%%%%%%%%%%%%%%%%%%%%%%
%
%%%%%%%%%%%%%%%%%%%%%%%%%%%%%%%%%%%%%%%%%%%%%%%%%%%%%%%%%%%%%%%%% Spray %%%%%%%%%%%%%%%%%%%%%%%%%%%%%%%%%%%%%%%%%%%%%%%%%%%%%%%%%%%%%%%
%
%%%%%%%%%%%%%%%%%%%%%%%%%%%%%%%%%%%%%%%%%%%%%%%%%%%%%%%%%%%%%%%%% Spray %%%%%%%%%%%%%%%%%%%%%%%%%%%%%%%%%%%%%%%%%%%%%%%%%%%%%%%%%%%%%%%
\section{Sprays}
%%%%%%%%%%%%%%%%%%%%%%%%%%%%%%%%%%%%%%%%%%%%%%%%%%%%%%%%%%%%%%%%%%
\begin{Def}
A spray $\zeta$ is a second order vector field on $M$ such that on
a local chart $(U,\varphi)$ it is determined by a smooth mapping
$\Gamma_\varphi:\varphi U\longrightarrow L^2_s(\mathbb{E},
\mathbb{E})$ in the following way:
\begin{eqnarray*}
[\zeta(v)]_\varphi(\varphi{m},v_\varphi)=(v_\varphi,
\Gamma_{\varphi}(\varphi m)(v_\varphi,v_\varphi)); ~m\in U, ~ v\in
T_mM
\end{eqnarray*}
(see \cite{Kumar}). Note that this definition coincides with the
one given in \cite{lang}.
\end{Def}

%%%%%%%%%%%%%%%%%%%%%%%%%%%%%%%%%%%%%%%%%%%%%%%%%%%%%%%%

\begin{The}
The limit of a  projective system of sprays on $M_i$ is a spray on
$M=\varprojlim M_i$.
\begin{proof}
For every $i \in \mathbb{N}$, let $\zeta_i$ be a second order vector
field on $M_i$. Moreover suppose that $(\varprojlim
U_i,\varprojlim \varphi_i)$ is a chart of $M=\varprojlim M_i$.
Then on the chart $(U_i,\varphi_i)$ on $M_i$, $\zeta_i$ is
determined by the  map
$\Gamma_{\varphi_i}:\varphi_i{U_i}\longrightarrow
L^2(\mathbb{E}_i,\mathbb{E}_i)$ with the property
\begin{eqnarray*}
[\zeta_i(v_i)]_{\varphi_i}({\varphi_i}m_i,v{_{\varphi_i}})=({v}_{\varphi_i},
{\Gamma}_{\varphi_i}(\varphi_i{m_i})({v}_{\varphi_i},{v}_{\varphi_i}));
 ~m_i\in U_i, ~v_i \in T_{m_i}M_i.
\end{eqnarray*}
To prove the result, it suffices to check that for $j\geq i$,
%%%%%%%%%%%%%%%%%%%%%%%%
$$(\rho_{ji}\times \rho_{ji})[\zeta_j(v_j)]_{\varphi_j}=
[\zeta_i(v_i)]_{\varphi_i}(\rho_{ji}\times \rho_{ji}).$$ Indeed
for every  $m_j \in U_j$ and $v_j=[U_j,\varphi_j,v_{\varphi_j}] \in
T_{m_j}{M_j}$ one obtains;
 %%%%%%%%%%%%%%%%%%%%%%%
\begin{eqnarray*}
(\rho_{ji}\times
\rho_{ji})[\zeta_j(v_j)]_{\varphi_j}(\varphi_j{m_j},v_{\varphi_j})
=%%%%%%%%%%%%%%%%%%%%%%%%%
(\rho_{ji}\times \rho_{ji})({v}_{\varphi_j},
{\Gamma}_{\varphi_j}(\varphi_j{m_j})({v}_{\varphi_j},{v}_{\varphi_j}))\\
%%%%%%%%%%%%%%%%%%%%%%%
&&\hspace{-11cm}=(v_{\varphi_i},\Gamma_{\varphi_i}(\varphi_i{m_i})(\rho_{ji}\times
\rho_{ji})(v_{\varphi_j},v_{\varphi_j}))
%%%%%%%%%%%%%%%%%%%%%
=({v}_{\varphi_i},
{\Gamma}_{\varphi_i}(\varphi_i{m_i})({v}_{\varphi_i},{v}_{\varphi_i}))\\
%%%%%%%%%%%%%%%%%%%%%%%%
&&\hspace{-11cm}=[\zeta_i(v_i)]_{\varphi_i}({\varphi_i}{m_i},v{_{\varphi_i}})
%%%%%%%%%%%%%%%%%%%%%%%
=[\zeta_i(v_i)]_{\varphi_i}(\rho_{ji}(\varphi_j{m_j}),\rho_{ji}({v}_{\varphi_j}))\\
&&\hspace{-11cm}=[\zeta_i(v_i)]_{\varphi_i}(\rho_{ji}\times
\rho_{ji})(\varphi_j{m_j},v_{\varphi_j}).
%%%%%%%%%%%%%%%%%%%%
\end{eqnarray*}
As mentioned in \cite{Kumar} if $\zeta_i$ is a spray on
$M_i$, for every pair of charts $(U_i,\varphi_i)$ and $(V_i,\psi_i)$
of $M_i$ at $m_i$, the transformation formula for
$\Gamma_{\varphi_i}$ is
\begin{eqnarray*}
\Gamma_{\psi_i}(\psi_im_i)(v_{\psi_i},v_{\psi_i})=D^2F_i(\varphi_i)
(v_{\varphi_i},v_{\varphi_i})+DF_i(\varphi_im_i).\Gamma_{\varphi_i}(\varphi_im_i)(v_{\varphi_i},v_{\varphi_i})
\end{eqnarray*}
where $F_i=\psi \circ \varphi^{-1}$ and
$v_i=[U_i,\varphi_i,v_{\varphi_i}] \in T_{m_i}M_i$. Suppose that
$\zeta=\varprojlim \zeta_i$ be a spray on $M=\varprojlim M_i$.
Then for charts $(U$=$\varprojlim U_i$, $\varphi$=$\varprojlim
\varphi_i)$ and $(V$=$\varprojlim V_i$, $\psi$=$\varprojlim
\psi_i)$ at $m=(m)_{i \in \mathbb{N}} \in M$ and
$v=[U,\varphi,v_{\varphi}] \in T_mM$:
\begin{eqnarray*}
\Gamma_{\psi}(\psi m)(v_{\psi},v_{\psi})&=&\varprojlim
\Gamma_{\psi_i}(\psi_im_i)(v_{\psi_i},v_{\psi_i})\\
&=&\varprojlim [D^2F_i(\varphi_im_i)
(v_{\varphi_i},v_{\varphi_i})+DF_i(\varphi_im_i).
\Gamma_{\varphi_i}(\varphi_im_i)(v_{\varphi_i},v_{\varphi_i})]\\
&=&D^2F(\varphi m) (v_{\varphi},v_{\varphi})+DF(\varphi m).
\Gamma_{\varphi }(\varphi m)(v_{\varphi},v_{\varphi})
\end{eqnarray*}
It means that the spray $\zeta=\varprojlim \zeta_i$ defines the
Christffel structure  $\{\Gamma_{\varphi}=\varprojlim
\Gamma_{\varphi_i}\}$ on $M=\varprojlim M_i$.
\end{proof}
\end{The}

\section{Dissections}
The concept of dissection is considered next. Kumar and Viswanath~\cite{Kumar}
established a one-to-one correspondence between dissections of $M$ and
Christoffel structures on $M$ for a Banach manifold $M$.
We extend this correspondence to projective
limit manifolds.

%%%%%%%%%%%%%%%%%%%%%%%%%%%%%%%%%%%%%%%%%%%%%%%%%%%%%%%%%

For $m \in M$, let $G_m:=\{f \in C^{\infty}(U_m)$ : $U_m$ is a
neighbourhood of $m \}$ and $G^{0}_m :=\{f \in G_m : f(m)=0\}$.
Define the space of 1-jets at $m$, denoted by $J_mM$, to be the
set of all equivalence classes in $G^{0}_m$, where two functions
$f,g \in G^{0}_m$ are equivalent if on every chart $(U, \varphi)$ of
$M$, the following relation holds true: $Df_\varphi(\varphi
m)=Dg_\varphi(\varphi m)$. In a similar way for every chart $(U,
\varphi)$ of $M$, one may define $J^2_mM:=\{[f] \in J_mM :
D^2f_\varphi(\varphi m)=D^2g_\varphi(\varphi m), \forall g \in [f]
\}$.
%%%%%%%%%%%%%%%%%%%%%%%%%%%%%%%%%%%%%%%%%%%%%%%%%%%
If $s \in J^2_mM$, then the local representation of $s$ on the
chart $(U,\varphi)$ is $s_{\varphi}=\alpha_{\varphi}\oplus
B_{\varphi} \in \mathbb{E}^*\oplus L^2_s(\mathbb{E},\mathbb{R})$
with transformation rule $\alpha_{\psi}=\alpha_{\varphi} \circ
DG(v)$ and $B_{\psi}=B_{\varphi}\circ (DG(v)\times
DG(v))+\alpha_{\varphi}\circ DG(v) \circ D^2F(u)\circ (DG(v)\times
DG(v))$, where $\alpha_{\varphi}$ is the local representation of
$\alpha \in T^*_mM$, $G=\varphi \circ \psi^{-1}$, $u=\varphi m$
and $v=\psi m$ (for more details see \cite{Kumar}).
%%%%%%%%%%%%%%%%%%%%%%%%%%%%%%%%%%%%%%%%%%%%%%%%%%%%%%%%%%
\begin{Def}    \label{a}
A dissection on $M$ is a map that to every $m \in M$ assigns a
closed subgroup of $J^2_{m}M$ say $D_m$. This is done in such a way that for
every chart $(U, \varphi)$ there exists a smooth mapping
$\Gamma_{\varphi}:\varphi U\longrightarrow
L^2_s(\mathbb{E},\mathbb{E})$ such that
$B_{\varphi}=\alpha_{\varphi}\circ \Gamma_{\varphi}(u)$ for $s \in
D_m$ and $s_{\varphi}=\alpha_{\varphi}\oplus B_{\varphi}$. In
other words
%\begin{eqnarray*}
$[D_m]_\varphi=\{\alpha\oplus{\alpha\circ{\Gamma_{\varphi}(u)}} :
\alpha\in \mathbb{E}^* \}$ (\cite{Kumar}).
%\end{eqnarray*}
\end{Def}
%%%%%%%%%%%%%%%%%%%%%%%%%%%%%%%%%%%%%%%%%%%%%%%%%%%%%%
We extend Kumar and Viswanath's results to projective limit
Fr\'{e}chet manifolds.
\begin{Pro}
If $\{M_{i}\}_{i \in \mathbb{N}}$ is a projective system of
manifolds and $\varprojlim J^2_{m_i}M_i$ exists then $\varprojlim
J^2_{m_i}M_i=J^2_{(m_i)}{\varprojlim M_i} $ (set-theoretically).
\begin{proof}
Let $G_m:=\{(f_i)_{i \in \mathbb{N}}; f_i:U_{m_i} \longrightarrow
\mathbb{R}$ is continuous  and $\varprojlim f_i$ exists$\}$ and
$G^0_m:=\{(f_i)_{i \in \mathbb{N}}\in G_m : f_i(m_i)=0,  ~\forall
i \in \mathbb{N}\}$. By defining
\begin{eqnarray*}
p : J^2_{m}M &\longrightarrow& \varprojlim
J^2_{m_i}M_i\\
{[f,m]} &\longmapsto& ([f_i, m_i]_i)_{i \in \mathbb{ N}}
\end{eqnarray*}
It can be checked that $p$ is well defined; moreover,  $p$ is one
to one since $p[f,m]=p[g,m]$ yields
\begin{eqnarray*}
{[f_i,m_i]}_i=[g_i,m_i]_i ,~ i \in \mathbb{N}.
\end{eqnarray*}
Hence
%\begin{eqnarray*}
$[f,m]=[\varprojlim f_i,(m_i)_{i \in \mathbb{N}}]_i=
\varprojlim[f_i,m_i]_i=\varprojlim [g_i,m_i]_i=[\varprojlim
g_i,(m_i)_{i \in \mathbb{N}}]=[g,m].$\\
%\end{eqnarray*}
Furthermore  $p$ is surjective. In fact if $([f_i,m_i]_i)_{i \in
\mathbb{N}}$ is an arbitrary element of $\varprojlim
{J^2_{m_i}M_i}$, we define $a=[\varprojlim f_i,(m_i)_{i \in
\mathbb{N}}]$. Then $p(a)=([f_i,m_i]_i)_{i \in \mathbb{N}}$ and
therefore $p$ is an isomorphism between $J^2_{m}M$ and
$\varprojlim J^2_{m_i}M_i$.
\end{proof}
\end{Pro}

%%%%%%%%%%%%%%%%%%%%%%%%%%%%%%%%%%%%%%%%%%%%%%%%%%%%%%%%%%%%%%%%
\begin{The}
The limit of a projective system of dissections of $\{M_i\}_{i \in
\mathbb{N}}$ is a dissection of $\varprojlim M_i=M_i$.
\end{The}
\begin{proof}
For every $i \in \mathbb{N}$, suppose $D_{m_i}$ is the closed
subgroup of $J^2_{m_i}{M_i}$ with the above mentioned  properties.
Moreover for $j \geq i$,
\begin{eqnarray*}
&&B_{\varphi_j}=\alpha_{\varphi_j}\circ \Gamma_{\varphi_j}({u_j})=
(\alpha_{\varphi_i} \circ \rho_{ji})\circ
\Gamma_{\varphi_j}({u_j})= \alpha_{\varphi_i} \circ
(\Gamma_{\varphi_i}(u_i)\circ (\rho_{ji}\times
\rho_{ji}))\\
&&=B_{\varphi_i}\circ (\rho_{ji}\times \rho_{ji}).
\end{eqnarray*}
Therefore $\varprojlim D_{m_i}$ exists and it is a dissection on
$M=\varprojlim M_i$.
\end{proof}
%%%%%%%%%%%%%%%%%%%%%%%%%%%%%%%%%%%%%%%%%%%%%%%%%%%%%%%%%%%
If $\varprojlim D_{m_i}$ is a dissection of $\varprojlim M_i=M$
and $(U$=$\varprojlim U_i$, $\varphi$=$\varprojlim \varphi_i)$,
$(V$=$\varprojlim V_i$, $\psi$=$\varprojlim \psi_i)$ are two
charts at $m=(m_i)_{i \in \mathbb{N}} \in M$, then
\begin{eqnarray*}
\Gamma_{\psi}(v)=\varprojlim \Gamma_{\psi_i}(v_i)=\varprojlim
[D^2F_i(u_i)\circ (DG_i(v_i)\times DG_i(v_i))\\+DF_i(u_i)\circ
\Gamma_{\varphi_i}(u_i) \circ (DG_i(v_i) \times DG_i(v_i))]\\
=D^2F(u)\circ (DG(v)\times DG(v))+DF(u)\circ \Gamma_{\varphi}(u)
\circ (DG(v) \times DG(v)),
\end{eqnarray*}
which precisely coincides with the Christoffel structures
$\{\Gamma_{\varphi}=\varprojlim \Gamma_{\varphi_i}\}$. (For more
details see \cite{Kumar}.) Hence we get the following result.

\begin{Cor}
There is one-to-one correspondence between dissections and
Christoffel structures on $M=\varprojlim M_i$.
\end{Cor}
%%%%%%%%%%%%%%%%%%%%%%%%%%%%%%%%%%%%%%%%%%%%%%%%%%%%%%%  Applications-Examples          %%%%%%%%%%%%%%%%%%%%%%%%%%%%%%%%%%%%%
%
%%%%%%%%%%%%%%%%%%%%%%%%%%%%%%%%%%%%%%%%%%%%%%%%%%%%%%%  Applications-Examples          %%%%%%%%%%%%%%%%%%%%%%%%%%%%%%%%%%%%%
%
%%%%%%%%%%%%%%%%%%%%%%%%%%%%%%%%%%%%%%%%%%%%%%%%%%%%%%%  Applications-Examples          %%%%%%%%%%%%%%%%%%%%%%%%%%%%%%%%%%%%%
\section{Examples}
%
%%%%%%%%%%%%%%%%%%%%%%%%%%%%%%%%%%%%%%%%%%%%%%%%%

\begin{Examp}
{\bf{The direct connection}}\\
 Let $G$ be a Banach Lie group with
the model space $\mathbb{E}$. Consider the mapping $\mu :
G\times{\mathbb{g}}\longrightarrow TG$ given by
$\mu(m,v)=T_e{\lambda_{m}}(v)$, where $\lambda_m$ is the left
translation on $G$ and $\mathbb{g}$ is the Lie algebra of $G$.
According to Vassiliou~\cite{Vass}, there exists a unique connection
$\nabla^{G}$ on $G$ which is $(\mu,id_G)-$related to the canonical
flat connection on the trivial bundle
$L=(G\times{\mathbb{g}},pr_1,G)$. Locally the Christoffel symbols
$\Gamma^G$ of $\nabla^G$ are given by
\begin{eqnarray*}
\Gamma^{G}_{\varphi}(x)(a,b)=Df_{\varphi}(x)(a,f^{-1}_{\varphi}(m)(b));~{x
\in \varphi{U}},  ~{a,b \in \mathbb{E}}
\end{eqnarray*}
where $f_\varphi$ is the local expression of the isomorphism
$T_e\lambda_{x}:T_eG\longrightarrow T_xG$ and $(U,\varphi)$ chart
of $G$. If $G=\varprojlim G_i$ is obtained as a projective limit of Banach
Lie groups and $\nabla^{G_i}$ is the direct connection on
$L^i=(G_{i}\times{\mathbb{g}}_i,pr_1,G_i)$, then
$\nabla^G=\varprojlim \nabla^{G_i}$ is exactly the direct
connection on $L=(\varprojlim G_{i}\times\varprojlim
{\mathbb{g}}_i,pr_1,\varprojlim G_i)$ \cite{PLBL}. Also, $\nabla^G$
determines a unique spray on $G=\varprojlim G_i$ locally given by
\begin{eqnarray*}
[\zeta^{G}(v)]_\varphi(\varphi{m},v_\varphi)=(v_\varphi,
{\Gamma}_{\varphi}^{G}(\varphi m)(v_\varphi,v_\varphi));~m\in U,
~v\in T_mG.
\end{eqnarray*}
Moreover, using $\nabla^G$ the Christoffel structure
$\{\Gamma_{\varphi}\}$ and  Hessian structure $H^G$ are obtained
where $H^G$ is locally given by
\begin{eqnarray*}
[H^{G}f(X,Y)]_\varphi(\varphi m)
%%%%%%%%%%%%%%%%%%%%
 =D^2f_\varphi(\varphi m)(X_\varphi
(\varphi m),Y_\varphi (\varphi m))+\\
%%%%%%%%%%%%%%%%%%%%
Df_\varphi(\varphi m).{\Gamma}_{\varphi}^{G}(\varphi m)(X_\varphi
(\varphi m),Y_\varphi (\varphi m)).
\end{eqnarray*}
\end{Examp}

%%%%%%%%%%%%%%%%%%%%%%%%%%%%%%%%%%%%%%%%%%%%%%%%%%%%

\begin{Examp}
{\bf{The flat connection}}\\
Let $M=\mathbb{E}$ with the global chart
$(\mathbb{E},id_{\mathbb{E}})$. The canonical flat connection
$\nabla^C$ on the trivial bundle $(M\times \mathbb{E},pr_1,M)$ is
locally  given by the Christoffel structure $\{\Gamma^C\}$, where
$\Gamma^C(x)(u)=0$, for every $(x,u) \in \mathbb{E}\times
\mathbb{E}$. Let $M=\mathbb{F}=\varprojlim \mathbb{E}_i$  and
consider it with the global chart
$(\mathbb{F},id_{\mathbb{F}})=\varprojlim (\mathbb{E}_i,
id_{\mathbb{E}_i})$. For the canonical flat
connection ${\Gamma^{C}}=\varprojlim \Gamma_i^C$ on $(M\times
\mathbb{F},pr_1,M),$ the spray $\zeta^C$ and the Hessian
structure $H^C$ are given by
\begin{eqnarray*}
[\zeta^{C}(v)]_\varphi(\varphi{m},v_\varphi)=(v_\varphi); ~m\in U,
~v\in T_mM
\end{eqnarray*}
and
\begin{eqnarray*}
[H^{C}f(X,Y)]_\varphi(\varphi m)
%%%%%%%%%%%%%%%%%%%%
 =D^2f_\varphi(\varphi m)(X_\varphi
(\varphi m),Y_\varphi (\varphi m)).
\end{eqnarray*}
\end{Examp}

%%%%%%%%%%%%%%%%%%%%%%%%%%%%%%%%%%%%%%%%%%%%%%%%%%%%%%%%%%%%

\section{Ordinary differential equations}
A curve $\gamma :(-\varepsilon,\varepsilon) \longrightarrow M$ is
called autoparallel or a geodesic with respect to the connection $\nabla$ if
$\nabla_{T{\gamma}}T{\gamma}=0$ (\cite{Vi}). Let $(U,\varphi)$ be
a local chart on $M$ and set
$\gamma_\varphi:=\varphi\circ\gamma:(-\varepsilon,\varepsilon)
\longrightarrow \mathbb{E} $, $\gamma^\prime_{\varphi}(t):=
T{\gamma_\varphi}:(-\varepsilon,\varepsilon) \longrightarrow
T\mathbb{E}$.

In this case the local  expression of
$\nabla_{T{\gamma}}T{\gamma}=0$ takes the form:
\begin{eqnarray*}
{\nabla_{T{\gamma_{\varphi}}}T{\gamma}}_{\varphi}(\gamma_{\varphi}(t)=
D\gamma^\prime_{\varphi}(t).\gamma^\prime_{\varphi}(t)-\Gamma_{\varphi}(\gamma_{\varphi}(t))
[\gamma^\prime_{\varphi}(t),\gamma^\prime_{\varphi}(t)]=0.
\end{eqnarray*}

Every spray is a second order vector field, hence every integral
curve of $\zeta$ is the canonical lifting of $\pi\circ\beta,$ so
$T(\pi\circ\beta)=\beta$. The curve
$\gamma:(-\varepsilon,\varepsilon)\longrightarrow M$ is a geodesic
spray with respect to $\zeta$ if $T\gamma$ is an integral curve
for $\zeta,$ namely, $T_{T_t\gamma(v_t)}T_t\gamma(v_t)=\zeta
T_t\gamma(v_t)$, where $v_t\in T_t{\mathbb{R}}$ with
$pr_2(v_t)=1$. In  local charts we have;
\begin{eqnarray*}
\big(\zeta(T_t\gamma{(v_t)})\big)_\varphi(\gamma_{\varphi}(t),{D_{t}}{\gamma_{\varphi}}(v_t))=
\big(\gamma_{\varphi}(t),\Gamma_{\varphi}(\gamma_{\varphi}(t))
[{D_{t}}{\gamma_{\varphi}}(v_t),{D_{t}}{\gamma_{\varphi}}(v_t)]\big).
\end{eqnarray*}
and
\begin{eqnarray*}
(T_{T_t\gamma(v_t)}T_t\gamma(v_t))_\varphi
=\big({D_{t}}{\gamma_{\varphi}}
(v_t),D_{{D_{t}}{\gamma_{\varphi}}(v_t)}D_t\gamma_{\varphi}(v_t,v_t)\big)
:=({\gamma}^{\prime}_{\varphi}(t),{\gamma}^{\prime\prime}_{\varphi}(t))
\end{eqnarray*}
So $\gamma$ must satisfy the (local)  equation
\begin{eqnarray*}
{\gamma}^{\prime\prime}_{\varphi}(t)=\Gamma_{\varphi}(\gamma_{\varphi}(t))
({\gamma}^{\prime}_{\varphi}(t),{\gamma}^{\prime}_{\varphi}(t)).
\end{eqnarray*}
Consequently the following theorem holds for Banach modelled
manifolds.
\begin{The}
Let $\zeta$ be the spray assigned to $\nabla$. There is a
one-to-one correspondence between geodesics of $\nabla$ and
geodesic sprays of $\zeta$.
\end{The}

Here we try to generalize this to the case of
Fr\'{e}chet manifolds where difficulties arise due to intrinsic
problems of the model spaces of these manifolds and mainly due to
the inability to solve general differential equations (see
\cite{Balkan}, \cite{Gal o} and \cite{Neeb}). We show that if
one focuses on the category of projective limit manifolds,
then similar results can be obtained.

%%%%%%%%%%%%%%%%%%%%%%%%%%%%%%%%%%%%%%%%%%%%%%%%%%%%%%%%%%%%%%%%%%%
\begin{The}
Let $M=\varprojlim M_i$ and $\zeta=\varprojlim \zeta_i$ be a spray
on $M$ with k-Lipschitz local components. Let $x_0\in M$ and
$y_0\in T_{x_0}M$. If for a chart $(U,\varphi)$ around $x_0$,
$M_\varphi=sup\{{\big{(}p_i(x_0)^2+p_i(\Gamma_\varphi(x_0)[y_0,y_0])^2\big{)}}^{1/2};~i\in
\mathbb{N}\}<\infty$, then there exists a locally unique geodesic
spray $\gamma :(-\varepsilon,\varepsilon)\longrightarrow M$ such
that $\gamma(0)=x_0$, $T_t\gamma(0)=y_0$ and $\varepsilon>0$ is
independent of the index $i$.
\begin{proof}
Let $\zeta:TM\longrightarrow TTM$ be a spray. Consider
$\{\zeta_i\}_{i \in \mathbb{N}}$, $x_0=({x_0}_{i})_{i \in
\mathbb{N}} \in \varprojlim M_i$ and $y_0=({y_0}_i)_{i \in
\mathbb{N}} \in \varprojlim T_{x_{0i}}M_i$. For every $i \in
\mathbb{N}$, $\zeta_i$ is a spray on $M_i$. Since $M_i$ is a
Banach manifold, by the existence theorem for ordinary differential
equations, there exists
$\gamma_i:(-\varepsilon_i,\varepsilon_i) \longrightarrow M_i$ with
\begin{eqnarray}
{\gamma_i}^{\prime\prime}_{\varphi_i}(t)=\Gamma_{\varphi_i}(\gamma_{\varphi_i}(t))
[{\gamma_i}^{\prime}_{\varphi_i}(t),{\gamma_i}^{\prime}_{\varphi_i}(t)],
\end{eqnarray}
satisfying  $\gamma_i(0)={x_0}_i$ and $T_{t\gamma_i}(0)={y_0}_i$.
For $j\geq i$, we claim that
$\varphi_{ji}\circ{\gamma_j}=\gamma_i$ and consequently
$\{\gamma_i\}_{i \in \mathbb{N}}$ forms a projective system of
curves on $\{M_i\}_{i \in \mathbb{N}}$ with the limit
$\gamma=\varprojlim \gamma_i$. Note that
\begin{eqnarray*}
({\varphi_i\circ\varphi_{ji}\circ\gamma_j}_{\varphi_j})^{\prime\prime}(t)
=(\rho_{ji}\circ{\varphi_j}\circ{\gamma_j}_{\varphi_j})^{\prime\prime}(t)
=\rho_{ji}({\varphi_j}\circ{\gamma_j}_{\varphi_j})^{\prime\prime}(t))
=\rho_{ji}\Gamma_{\varphi_j}({\gamma_j}_{\varphi_j}(t))\\
({\gamma_j}_{\varphi_j}^\prime(t),{\gamma_j}_{\varphi_j}^\prime(t))
=\Gamma_{\varphi_i}((\rho_{ji}\circ{\varphi_j}\circ{\gamma_j}_{\varphi_j})(t))
[(\rho_{ji}\circ{\varphi_j}\circ{\gamma_j}_{\varphi_j})^\prime(t),
(\rho_{ji}\circ{\varphi_j}\circ{\gamma_j}_{\varphi_j})^\prime(t)]\\
=\Gamma_{\varphi_i}(({\varphi_i\circ\varphi_{ji}\circ\gamma_j}_{\varphi_j})(t))
[({\varphi_i\circ\varphi_{ji}\circ\gamma_j}_{\varphi_j})^\prime(t),
({\varphi_i\circ\varphi_{ji}\circ\gamma_j}_{\varphi_j})^\prime(t)].
\end{eqnarray*}
Moreover
$(\varphi_{ji}\circ{\gamma_j})(0)=\varphi_{ji}({x_0}_j)={x_0}_i$
and $T_t(\varphi_{ji}\circ{\gamma_j})(0)={y_0}_i$. By uniqueness
of solutions for ordinary differential equations on Banach spaces
(manifolds) we have $\varphi_{ji}\circ{\gamma_j}={\gamma_i}$ and
consequently $\gamma=\varprojlim \gamma_i$ exists. Furthermore
\begin{eqnarray*}
T_{T_t\gamma(v_t)}T_t
\gamma(v_t)=\{T_{T_t\gamma_i(v_t)}T_t\gamma_i(v_t)\}_{i \in
\mathbb{N}}=\{{\zeta_i}(T_t\gamma_i(v_t))\}_{i \in
\mathbb{N}}=\zeta (T_t\gamma(v_t)).
\end{eqnarray*}
According to Theorem~\ref{ExandUniq}, $\varepsilon$ does not converge to $0$
and consequently there exists $\varepsilon>0$ such that $\gamma
:(-\varepsilon,\varepsilon)\longrightarrow M$ is a local geodesic
spray with respect to $\zeta$.

 Let
$\beta:(-\varepsilon_1,\varepsilon_1)\longrightarrow M$ be another
curve such that $T_{T_t\beta(v_t)}T_t\beta(v_t)=\zeta
(T_t\beta(v_t))$, satisfying in the above boundary conditions. For
every $i \in \mathbb{N}$, $\beta_i=\psi_i\circ\beta$ satisfies in
equation (2) with $\beta_i(0)={x_0}_i$ and
$T_{t}{\beta_i}(0)={y_0}_i$. Hence $\beta_i=\gamma_i$ and
consequently $\beta=\varprojlim \beta_i=\varprojlim
\gamma_i=\gamma$ on the intersection of their domains.
\end{proof}
\end{The}

Finally in a similar way one can prove the theorem for geodesics with
respect to the connection $\nabla$. As a conclusion we can state
the following corollary.
\begin{Cor}
For a projective limit manifold $M=\varprojlim M_i$ there is a
one-to-one correspondence between (linear) connections and sprays.
Moreover, the geodesics of $\nabla$ are geodesic sprays of
$\zeta$.
\end{Cor}

%%%%%%%%%%%%%%%%%%%%%%%%%%%%%%%%%%%%%%%%%%%%%%%%%%%%%%%%%%%%%%%%%%
\section{Parallel translation}
Vilms~\cite{Vi} defines a connection on $M$ as a vector bundle
morphism $\nabla :T(TM)\longrightarrow TM$. So $\nabla$ is fully
determined by its local components, called  Christoffel symbols,
denoted by $\{\Gamma_{\alpha}\}_{\alpha \in I}$ corresponding to
an atlas of charts $\{(U_\alpha,\varphi_\alpha)\}_{\alpha\in I}$
of $M$. Then, $\Gamma_{\varphi}:\varphi U\longrightarrow
L^2(\mathbb{E},\mathbb{E})$, and for two charts $(U,\varphi)$ and
$(V,\psi)$ at $m \in M$, $e_1,e_2 \in \mathbb{E}$,
 $u=\varphi(m), v=\psi(m)$, $F=\psi\circ {\varphi}^{-1}$ we have
\begin{eqnarray*}
\Gamma_{\psi}(v)((DF(u).e_1,DF(u).e_2)=DF(u).\Gamma_{\varphi}(u)(e_1,e_2)+D^2F(e_1,e_2).
\end{eqnarray*}
Clearly, our definition coincides with the above; we next consider
parallel transport of vectors along a curve.
%%%%%%%%%%%%%%%%%%%%%%%%%%%%%%%%%%%%%%%%%%%%%%%%%%%%%%%%%%%%%%
\begin{The}
Given $\nabla :T(TM)\longrightarrow TM$
 a connection on $(TM,M,\pi),$  take a smooth curve
$c:(a,b)\longrightarrow M$ with $0 \in (a,b),
 c(0)=x$. Then, there is a neighbourhood U of $T_xM\times\{0\}
\subseteq  T_xM\times (a,b)$ and a smooth mapping
%%%%%%%%%%%%%%%%%%%%%
$\bar{c}:U\longrightarrow TM$ such that:\\
(i) $\pi(\bar{c}(u_x,t))=c(t)$ and $\bar{c}(u_x,0)=u_x$,\\
(ii) $\nabla(\frac d{dt} (\bar{c})(u_x,t))=0.$
%%%%%%%%%%%%%%%%%%%%
\begin{proof}
For every $(U,\varphi)$ chart of $M$,   $\nabla(\frac
d{dt}\bar{c}(u_x,t))=0$, locally gives\\
$-\Gamma_{\alpha}(c(t))\big(\frac d{dt}c(t),\gamma(y,t)\big)+\frac
d{dt}\gamma(y,t)=0$, where
$T \varphi(\bar{c}(c,T{\varphi}^{-1}(x,y),t)):=(c(t),\gamma(y,t))$\\
(i.e. $\gamma:\mathbb{E}\times(a,b)\longrightarrow S)$. For $M$
a Banach manifold, by the existence theorem for differential
equations, $\bar{c}$ always exists.
\end{proof}
\end{The}

Using our method one can prove a similar theorem for
parallel transport along curves in the category of projective limit
manifolds. The equivalence of linear connections with sprays means
that parallel transport is equivalently determined by a spray \cite{Michor}.
%Interaction of such curves with Hessian structures.

\begin{Examp}
{\bf{Geodesics on the diffeomorphism group of the circle}}

The main reference for this example is Constantin and Kolev~\cite{Const-Kol}. Let
$D=Diff(\mathbb{S}^1)^+$ be the group of all smooth orientation-preserving
diffeomorphisms of the circle $\mathbb{S}^1$. We can endow $D$ with a smooth
manifold structure based on the
Fr\'{e}chet space $C^\infty(\mathbb{S}^1)$.\\
Moreover  a right invariant weak Riemannian metric on $D$ is
defined. Note that $C^\infty(\mathbb{S}^1)=\bigcap_{n\geq 2k+1}H^n(\mathbb{S}^1)$ where
$H^n(\mathbb{S}^1)$, $n\geq 0$ is the space $L^2(\mathbb{S}^1)$ of all square
integrable functions $f$ with the distributional derivatives up
to order $n$, $\partial_x^i$ with $i=0,1,...,n$, in $L^2(\mathbb{S}^1)$.
$H^n(\mathbb{S}^1)$, $n\geq 0$ is a Hilbert space with the norm
\begin{equation*}
\|f\|_n^2=\sum_{i=0}^n\int_S(\partial_x^if)^2(x)dx.
\end{equation*}
The main difference of this example for our method lies in the
existence of a metric and this allows us to prove the length
minimizing property of geodesics.

We move this problem to the projective limit framework. In this
special case like \cite{Omori1978} the connecting morphisms of the
model space are inclusions
\begin{eqnarray*}
{\rho}_{n+1,n}:H^{n+1}&\hookrightarrow & H^{n}\\
f&\longmapsto& f
\end{eqnarray*}
The meaning of this morphism is clear, namely if $f \in H^{n+1}$ with
the norm on $H^{n+1}$ then
$\sum_{i=0}^{n+1}\int_S(\partial_x^if)^2(x)dx < \infty$. Clearly
$\sum_{i=0}^n\int_S(\partial_x^if)^2(x)dx < \infty$, so the
function $f$ belongs precisely to $H^n$. Consequently if $f \in
\bigcap_{n\geq 2k+1}H^n(\mathbb{S}^1)$ then $(f)\in \varprojlim H^n(\mathbb{S}^1)$ and,
conversely, $C^\infty(\mathbb{S}^1)=\bigcap_{n\geq
2k+1}H^n(\mathbb{S}^1)=\varprojlim_{n\geq 2k+1} H^n(\mathbb{S}^1).$

For $\varphi_{0} \in D$ define $U_{0}=\{\varphi \in D :
{\parallel\varphi - \varphi_{0} \parallel}_{C^{0}(\mathbb{S}^1)} < 1/2 \}$
and $u:u_0\longrightarrow C^{\infty}(\mathbb{S}^1)$ with $u(x)=\frac 1{2\pi
{i}}ln \big(\overline{\varphi_{0}(x)}\varphi(x)\big)$, $x \in S$.
Then $(U_0,\psi_0)$ is a local chart with $\psi_0(\varphi)=u$ and
change of charts given by $\psi_2 \circ {\psi}^{-1}_1=u_1+ \frac
1{2\pi{i}}ln(\overline {\varphi_2}\varphi_1)$. Note that $\psi_2
\circ {\psi}^{-1}_1 : \psi_{1}(u_1)\subseteq C^\infty
(\mathbb{S}^1)\longrightarrow \psi_2(u_2)\subseteq C^\infty (\mathbb{S}^1)$  can be
recognized as the projective limit on Hilbert components, say
$(\psi_2 \circ {\psi}^{-1}_1)_i : H^i(\mathbb{S}^1)\longrightarrow H^i(\mathbb{S}^1)$ ,
$(\psi_2 \circ {\psi}^{-1}_1)=\varprojlim (\psi_2 \circ
{\psi}^{-1}_1)_i$. These maps are called k-Lipschitz and so
$(\psi_2 \circ {\psi}^{-1}_1)$. This structure endows $D$ with a
smooth manifold structure based on the Fr\'{e}chet space $C^
\infty (\mathbb{S}^1)$.

Let $k\geq 0$, for $n\geq 0$ define the linear seminorms $A_k:
H^{n+2k}(\mathbb{S}^1)\longrightarrow H^n(\mathbb{S}^1)$ with $A_k=1-{\frac {d^2}{dx^2}}
+...+(-1)^k \frac {d^{2k}}{dx^{2k}}.$ This enables us to define
the bilinear operator $B_k:C^{\infty} \times
C^{\infty}\longrightarrow C^{\infty}$ with
$B_k(u,v)={A_k}^{-1}(2v_xA_k(u)+vA_k(u_x)) ~ u,v \in C^{\infty}$.
Note that $B={\varprojlim}_{n\geq 2k+1}{B_k}^n$ where
${B_k}^n:H^n(\mathbb{S}^1) \times H^n(\mathbb{S}^1)\longrightarrow H^{n-2k}(\mathbb{S}^1).$
As stated in \cite{Const-Kol}, Theorem 1, there exists a unique
linear (Riemannian) connection $\nabla^{k}$ on $D$.

If $\varphi:J\longrightarrow D$ is a $C^2$-curve satisfying the
autoparallel equation with respect to the linear connection
$\nabla^k$, then
\begin{eqnarray*}
u_t=B_k(u,u), ~ t \in J
\end{eqnarray*}
where $u=\varphi_t \circ \varphi^{-1} \in T_{Id}D\equiv
C^{\infty}(\mathbb{S}^1)$. The term autoparallel rather than
geodesic is better since there is no underlying  Riemannian  metric. However
the utilization of a weak Riemannian metric is an issue that
remains open. Since $B_k=\varprojlim {B_k}^i,$ is the projective
system of bilinear maps (as Christoffel symbols) we can endow $D$
with the linear connection $\nabla_k=\varprojlim {\nabla_k}^i$. Given
an initial value we obtain the
unique autoparallel $\varphi :J\longrightarrow D$ obtained as the
projective limit on Hilbert components.  The problem is much
easier than the general case. Specifically, let the solution on
the $H^n(\mathbb{S}^1)$, $n\geq 2k+1$ have the manifold domain $[0,T_n)$ with
$T_n>0$. If $T_n\leq T_{\leq 2k+1}$ then $T_n=T_{2k+1}$ for all
$n\geq 2k+1$ i.e. the solution $\varphi_n$ on $H^n(\mathbb{S}^1)$, $n\geq 2k$
defined on $[0,T_{2k+1})$ for every $n$.

Note that in the general case there is no way to model the
diffeomorphism group of a manifold $M$ on a Banach space. However,
there is the possibility to view $Diff(M)$ as a projective limit
of Hilbert manifolds (\cite{Omori1974}). Moreover, the existence
of the geodesics in the general case of $Diff(M)$ is an open
question, so using the proposed technique with an
appropriate choice of imposed metric
may yield some results.
\end{Examp}

\section{Appendix: Existence and uniqueness theorem for second order ordinary differential equations on Fr\'{e}chet spaces}
Start with the assumptions of \cite{Gal-Pal}. Namely, let
$\mathbb{F}$ be a Fr\'{e}chet space and $\{p_i\}_{i\in\mathbb{N}}$
be a countable family of seminorms which determine the topology of
$\mathbb{F}$.

\begin{The}\label{ExandUniq}
Let $\mathbb{F}$ be a Fr\'{e}chet space and
$\Phi:\mathbb{R}\times\mathbb{F}\times\mathbb{F}\longrightarrow\mathbb{F}$
a projective limit $k$-Lipschitz mapping. For the second order
differential equation
\begin{equation}
x''=\Phi(t,x,x')
\end{equation}
with the initial condition $(t_0,x_0,y_0)$,  if there exists a
constant $\tau\in \mathbb{R}^+$ such that
\begin{equation*}
M=sup\{ {\big(p_i(y_0)^2+{p_i(\Phi(t,x_0,y_0))^2}\big)}^{1/2};~
i\in \mathbb{N}, t\in[t_0-\tau,t_0+\tau]\}<\infty
\end{equation*}
and $a=min\{\tau,{\frac{1}{M+k}}\}$, then (2) has a unique
solution on $I=[t_0-a,t_0+a]$.
\end{The}
\begin{proof}
If we set $x'=y$, $x'=y,~ y'=\Phi(t,x,y)$.
Denoting $z=(x,y)$ one takes:
\begin{equation}
z'=(x, y)^\prime=(x^\prime, y^\prime)=(y, \Phi(t,x,y)
)=\tilde{\Phi}(t,z)
\end{equation}
where
$\tilde\Phi:\mathbb{R}\times\mathbb{F}\times\mathbb{F}\times\mathbb{F}\longrightarrow\mathbb{F}\times\mathbb{F}$,
is  also a $k$-Lipschitz mapping. Since
\begin{eqnarray*}
{\big(p_i(y_0)^2+p_i{(\Phi(t,x_0,y_0))}^2\big)}^{1/2}=p_i(y_0,\Phi(t,x_0,y_0))=p_i(\tilde{\Phi}(t,z_0));
\end{eqnarray*}
and
\begin{equation*}
M=sup\{  p_i(\tilde{\Phi}(t,z_0))={\big(
{{p_i(y_0)}^2+p_i{(\Phi(t,x_0,y_0))}^2}\big)}^{1/2};~ i\in
\mathbb{N}, t\in[t-0-\tau,t-0+\tau]\}<\infty
\end{equation*}
 by Theorem 3 in \cite{Gal-Pal}, $(4)$ has a unique solution on
$I=[t_0-a,t_0+a]$ such that $a=min\{\tau,{1\over{M+k}}\}$. Hence
there exists also a solution for (3) say  $z:I\longrightarrow
\mathbb{F}\times\mathbb{F}$. If $z=(z_1,z_2)$ then, $z_1$ and
$z_2$ are unique solution for $x'=y$  and  $y'=\Phi(t,x,y)$
respectively on $I$. Consequently ${z_1}^{\prime}=y$,
$y'=\Phi(t,z_1,y)$ i.e.
\begin{equation*}
{z_1}''=\Phi(t,z_1,z_1')~on~ I.
\end{equation*}
Note that the interval $I$ is independent of the index $i$. For
each $i\in \mathbb{N}$ from the equation
\begin{equation*}
x_i''=\Phi_i(t,x_i,x_i')
\end{equation*}
with the initial condition $(t_0,{x_0}_i,{y_0}_i')$ we have the
unique solution $x_i$. On the other hand for $i\leq j$,
$f_{ji}\circ x_j$ is also a solution of (4) with $f_{ji}\circ
x_j(t_0)={x_0}_i$ and $(f_{ji}\circ x_j)'(t_0)={y_0}_i$. Hence
$f_{ji}\circ x_j=x_i$ for $i\leq j$, i.e. $x=\varprojlim x_i$ can
be defined. Moreover
\begin{equation*}
x''={(x_i'')}_{i\in\mathbb{N}}={(\Phi_i(t,x_i,x_i'))}_{i\in\mathbb{N}}=\Phi(t,x,x'),
\end{equation*}
i.e. $\varprojlim x_i$ is a solution for (2). The uniqueness of
$x$ follows from the uniqueness of solution for Banach components.
\end{proof}

%%%%%%%%%%%%%%%%%%%%%%%%%%%%%%%%%%%%%%%%%%%%%%%%%%%%%%%%%%%%%%%%%%%%%%%%%%%%%%%%%%%%%%%%%%%%%%%%%%%%%%%%%%%%%
%
%

\end{document}